\documentclass[11pt]{amsart}
\usepackage{mathdots, url}
\usepackage{graphicx}
\usepackage{amsmath, amsfonts, amssymb, amscd, amsthm}
\usepackage{yhmath}
\usepackage{setspace}
\usepackage{mathtools}
\usepackage[all]{xy}
\usepackage{verbatim}
\usepackage{enumerate}
\usepackage[usenames]{color}
\usepackage{varwidth}
\allowdisplaybreaks
\newtheorem{theorem}{Theorem}[section]

\newtheorem{lemma}{Lemma}[section]

\newtheorem{proposition}{Proposition}[section]

\theoremstyle{definition}
\newtheorem{definition}{Definition}[section]

\theoremstyle{remark}
\newtheorem{remark}{Remark}[section]

\numberwithin{equation}{section}
\allowdisplaybreaks
\begin{document}
\title[A problem on Hecke algebras for $\mathrm{GL}_n(F)$ for $n>2$ over $p$-adic field $F$]{A problem on Hecke algebras for $\mathrm{GL}_n(F)$ for $n>2$ over $p$-adic field $F$}
\author{Subha Sandeep Repaka}
\address{Department of Mathematics, SRM University- AP, Mangalagiri-Mandal, Neeru Konda, Amaravati, Andhra Pradesh 522240}
\email{subhasandeep.r@srmap.edu.in, sandeep.repaka@gmail.com}

	\date{\today}	
	\setcounter{tocdepth}{1}
	\date{\today}

	\begin{abstract}
		We study the Hecke algebra $\mathcal{H}_G(F)$ for $G = \mathrm{GL}_n$ and $n>2$ where $F$ is a non-Archimedean local field of characteristic zero. We show that for $G = \mathrm{GL}_n$ and $n>2$ and any two such fields $E$ and $F$, there is a Morita equivalence $\mathcal{H}_G(E) \sim \mathcal{H}_G(F)$, by using the Bernstein decomposition of the Hecke algebra and by determining the intertwining algebras that yield the Bernstein blocks up to Morita equivalence.
	\end{abstract}
	
		\maketitle
		
	\section{Introduction}
	
	Let $F$ be a non-Archimedean local field of characteristic 0. Define $G = \mathrm{GL}_n(F)$ where $n>2$. Let $C_c^\infty(G)$ be the space of locally constant compactly supported functions with values in $\mathbb{C}$. Let $\Phi_1, \Phi_2 \in \mathcal{H}_G(F)$ and let us choose a unimodular Haar measure $dg$ on $G$. We then define convolution given by $*$ on $C_c^\infty(G)$ via:
	
	\begin{equation}
		(\Phi_1 * \Phi_2)(h) = \int_G \Phi_1(hg^{-1}) \Phi_2(g)\, dg  \quad for \quad  h \in G
	\end{equation} 
	
   This gives the Hecke algebra $\mathcal{H}(G)$. Thus we have $\mathcal{H}(G)= C_c^\infty(G)$. We want to understand in this paper to what extend does the representation theory of $\mathcal{H}_G(F)$ determine the field $F$. We will show that for any two $p$-adic fields $E$ and $F$, there is a Morita equivalence $\mathcal{H}_G(E) \sim \mathcal{H}_G(F)$. Thus we prove the following theorem:
   
   \begin{theorem}\label{The_0}
   	Let $F,E$ be a non-Archimedean fields of characteristic 0. Let $G= \mathrm{GL}_n(F)$ where $n>2$. Let $\mathcal{H}_G(E)$ and $\mathcal{H}_G(F)$ and be local Hecke algebras over fields $E$ and $F$ respectively. Then Morita equivalence $\mathcal{H}_G(E) \sim \mathcal{H}_G(F)$ would imply $E \cong F$.
	\end{theorem}
	
	The proof relies on the Bernstein decomposition of the category of smooth representations together with the theory of types and covers. Using the description of Bernstein blocks in terms of affine Hecke algebras, we compare
	the module categories associated to different local fields and analyze the dependence of the Hecke algebra on the residue field cardinality and torsion parameters.
	
	The paper is organized as follows. In Section 2, we review Hecke algebras and smooth representations of \(p\)-adic groups. Section 3 discusses affine Hecke algebras and their defining relations. In Section 4, we recall
	parabolic induction and supercuspidal representations. Section 5 is devoted
	to the Bernstein decomposition of the category of smooth representations. In Sections 6 and 7, we study Bernstein blocks of local Hecke algebras, types, spherical Hecke algebras, and their relation with affine Hecke algebras, culminating in the proof of the main theorem $\ref{The_0}$. In section 8, main proof of the theorem $\ref{The_0}$ is discussed.
	
	\section{Hecke Algebra and Representations}
	
	If $(\pi, V)$ is a smooth complex representation of $G$  and $\Phi \in \mathcal{H}_G(F)$, we define the operator $\pi(\Phi) : V \longrightarrow V$ via:
	\begin{equation}
		\pi(\Phi)v = \int_G \Phi(g)\pi(g)v\, dg \quad for \quad  v \in V
	\end{equation}
	
	The above operator satisfies:
	 \begin{equation}
	 	\pi(\Phi_1)\pi(\Phi_2)= \pi(\Phi_1*\Phi_2) \quad for \quad  \Phi_1,\Phi_2 \in \mathcal{H}_G(F)
	 \end{equation}
	 
	Therefore, $V$ becomes an $\mathcal{H}_G(F)$-module and we can show that the smoothness of the representation $\pi$ makes $V$ a non-degenerate $\mathcal{H}_G(F)$-module. By saying $V$ is a non-degenerate module, we mean $\mathcal{H}_G(F)V=V$. Conversely, we can associate  with any non-degenerate
	$\mathcal{H}_G(F)$-module a unique smooth complex representation of $G$. Let us denote $\mathfrak{R}(G)$ to be the category of smooth complex representations of $G$ and let $\mathcal{H}_G(F)$-Mod denote the category of non-degenerate $\mathcal{H}_G(F)$-modules. Thus we have the following equivalence of categories:
	\begin{equation} \label{eqn_0}
		\mathfrak{R}(G) \simeq \mathcal{H}_G(F)\text{-Mod}.
	\end{equation}
	
	\section{Affine Hecke Algebra}
	
	Observe that the affine Hecke algebra $\mathcal{H}(r,z)$ for $r \in \mathbb{Z_{+}}$ and $z \in \mathbb{C}^{\times}$ is generated by the generators $s_i (1 \leq i \leq n), t, t^{-1}$ and is defined by the following relations: (cf.[\cite{BushnellKutzko1993}, Definition 5.4.6])
	
		\begin{enumerate}[(R1)]
			\item $ (s_i + 1)(s_i - r) = 0 \quad \text{for } 1 \leq i \leq m-1.$
			\item $t^2 s_1 = s_{m-1} t^2.$
			\item 	$ts_i = s_{i-1} t \quad \text{for } 2 \leq i \leq m-1.$
			\item	$ s_i s_{i+1} s_i = s_{i+1} s_i s_{i+1} 
			\quad \text{for } 1 \leq i \leq m-2.$
     	  	\item $s_i s_j = s_j s_i 
			\quad \text{for } 1 \leq i,j \leq m-1 \text{ with } |i-j| \geq 2.$
		\end{enumerate}
	\section{Parabolic Induction}
	\begin{definition}
		A quasicharacter $\chi$ of $F^\times$ is a continuous homomorphism
		\[
		\chi : F^\times \to \mathbb{C}^\times.
		\]
		It is called \emph{unramified} if it is trivial on $\mathcal{O}_F^\times$. Any unramified quasicharacter is of the form $|\cdot|^z$ for some $z \in \mathbb{C}$.
	\end{definition}
	
	\begin{lemma}(cf. [\cite{Wedhorn2008}, (2.1.18)])
		Every irreducible admissible representation $\pi$ which is finite-dimensional is one-dimensional and there exists a quasicharacter $\chi$ such that
		\[
		\pi(g) = \chi(\det g) \quad \text{for all } g \in G.
		\]
	\end{lemma}
	
	\begin{definition}
		A parabolic subgroup $P$ of $G$ is such that $G/P$ is complete. Equivalently, $P$ contains a Borel subgroup $B$. Every $P$ is the semidirect product of its unipotent radical and a $F$-closed reductive group $L$, called the Levi subgroup.
	\end{definition}
	
	\begin{remark}
		The proper parabolic subgroups of $\mathrm{GL}_n(F)$ are the block upper triangular matrices and their conjugates. For $n=2$, these are precisely the Borel subgroups.
	\end{remark}
	
	\begin{definition}
		Let $\rho$ be a smooth representation of a Levi subgroup $L$ of a parabolic subgroup $P$ of $G$. The parabolic induction $\iota_P^G(\rho)$ consists of locally constant functions $\phi$ on $G$ satisfying:
		
		\[
		\phi(pg) = \delta_P(p)^{1/2} \rho(p)\phi(g), \quad \text{for} \quad p \in P, g \in G.
		\]
		
		 The normalising factor $\delta_P = \triangle_P^{-1}$ is the inverse of the modular character $\triangle_P^{-1}$. The modular character satisfies:
		 
		 \[\triangle_P (diag(a_1, . . . , a_n)) =\mid a_1 \mid^{1-n}\mid a_2 \mid^{3-n} . . . \mid a_n \mid^{n-1}\]
		 
		 Note that parabolic induction preserves smoothness and admissibility but not necessarily irreducibility.
	\end{definition}
	
	The following two definitions are taken from \cite{Karemaker2016}:
	
	\begin{definition}\label{def_1}
		An infinite-dimensional irreducible admissible representation 
		$\pi : G \to \mathrm{GL}(V)$ is called (absolutely) cuspidal or 
		supercuspidal if it is not a subquotient of a representation that is 
		parabolically induced from a proper parabolic subgroup of $G$.
	\end{definition}

	\begin{definition}\label{def_2}
		A partition $(n_1,\dots,n_r)$ of $n$ means a partition of $\{1,2,\dots,n\}$ into segments $(1,\dots,n_1),\ (n_1+1,\dots,n_1+n_2),\ \dots,\ (n_1+\cdots+n_{r-1}+1,\dots,n)$ of respective lengths $n_i$. We will write $(n_1,\dots,n_r)\perp n$ for such a partition. For any $n_i$ appearing in a partition of $n$, write\[\Delta_i = \{\rho_i,\ \rho_i|\cdot|,\ \dots,\ \rho_i|\cdot|^{n_i-1}\}\]for $i=1,\dots,r$, where $\rho_i$ is an irreducible supercuspidal representation of $\mathrm{GL}_{n_i}(F)$.The $\Delta_i$ are also called segments, and we say that $\Delta_i$ precedes $\Delta_j$ if $\Delta_i \nsubseteq \Delta_j$ and $\Delta_j \nsubseteq \Delta_i$, if $\Delta_i \cup \Delta_j$ is also a segment, and $\rho_i = \rho_j |\cdot|^k$ for some $k>0$.
		\end{definition}
		
	We have the following Lemma based on (cf. [\cite{Zelevinsky1980}, Theorem6.1],
	[\cite{BernsteinZelevinsky1976}, Corollary 3.27], and [\cite{PrasadRaghuram2008}, 189–190 pp.]):
		\begin{lemma}
			For any partition $(n_1,\dots,n_r)$ of $n$ and a choice of segments
			so that $\Delta_i$ and $\Delta_{i+1}$ $(i=1,\dots,r)$ do not precede each other,
			there exists a corresponding induced representation, denoted
			\[
			\iota_P^G(\rho_1 \otimes \cdots \otimes \rho_r),
			\]
			whose unique irreducible quotient is an irreducible admissible representation of $G$.
			Any irreducible admissible representation of $G$ is equivalent to such a quotient representation.
		\end{lemma}
		
	From now on we denote $\rho_1 \otimes \cdots \otimes \rho_r$ by $\rho$.\par
	
	\section{Bernstein Decomposition}
	Recall that $G=\mathrm{GL}_n(F)$. According to Theorem 3.3 in \cite{Kutzko1998}, we have the following Proposition. 
	
	\begin{proposition} \label{cuspidal_pair}
		\begin{enumerate}
			\item Let $L$ be a Levi subgroup of $G$ (i.e., a Levi component of a parabolic subgroup $P$ of $G$). Let $\sigma$ be an irreducible smooth supercuspidal representation of $L$. Then $\iota_P^G{\sigma}$ has finite length for every parabolic subgroup $P$ with Levi component $L$. Further, the set of the composition factors or irreducible sub-quotients of $\iota_P^G\sigma$ is independent of $P$. 
			
			\item Let $L_1, L_2$ be Levi subgroups of $G$ and $\sigma_1, \sigma_2$ be irreducible supercuspidal  smooth representations of $L_1, L_2$ respectively. Then for any parabolic subgroups $P_1, P_2$ with Levi components $L_1, L_2$ respectively, we have the representations $\iota_{P_1}^G\sigma_1, \iota_{P_2}^G\sigma_2$ either have the same set of composition factors or have no composition factors in common. Now the representations  $\iota_{P_1}^G\sigma_1$ and $ \iota_{P_2}^G\sigma_2$ have the same set of composition factors $\Longleftrightarrow$ the pairs $(L_1, \sigma_1)$ and $(L_2, \sigma_2)$ are conjugate; that is, there is an element $g\in G$ such that $L_2=L_1^g= g^{-1}L_1g$ and $\sigma_2 \simeq {\sigma_1}^g$.
			
			\item Let $(\pi, V)$ be an irreducible smooth representation of $G$. Then there exists a parabolic subgroup $P$ of $G$ with Levi component $L$, unipotent radical $U$ and an irreducible supercuspidal smooth representation $\sigma$ of $L$ such that $\pi$ is equivalent to an irreducible sub-quotient or a composition factor of $\iota_P^G\sigma$. We refer to the pair $(L, \sigma)$ where $L$ is a Levi subgroup of $G$ and $\sigma$ is an irreducible supercuspidal smooth representation of $L$ as a cuspidal pair.
		\end{enumerate}
	\end{proposition}\par
	
	Now by Propn.~\ref{cuspidal_pair}, there exists unique conjugacy class of cuspidal pairs $(L, \sigma)$ with the property that $\pi$ is isomorphic to a composition factor of $\iota_P^G\sigma$ for some parabolic subgroup $P$ of $G$. We call this conjugacy class of cuspidal pairs, the cuspidal support of $(\pi, V)$.\par
	
	Given two cuspidal supports $(L_1, \sigma_1)$ and $(L_2, \sigma_2)$ of $(\pi, V)$, we say they are inertially equivalent if there exists $g \in G$ and $\chi \in \mathrm{X}_{nr}(L_2)$ such that $L_2= L_1^g$ and $\sigma_1^g\simeq \sigma_2\otimes\chi$. We write $[L, \sigma]_G$ for the inertial equivalence class or inertial support of $(\pi, V)$. Let $\mathfrak{B}(G)$ denote the set of inertial equivalence classes $[L, \sigma]_G$.\par
	
	We say that a smooth irreducible representation $(\pi, V)$ has inertial support $[L, \sigma]_G$ if $(\pi, V)$ appears as a sub-quotient of a representation parabolically induced from some element of $[L, \sigma]_G$. \par 
	
	Let $\mathfrak{R}(G)$ denote the category of smooth representations of $G$. Let ${\mathfrak{R}}^s(G)$ be the full sub-category of smooth representations of $G$ with the property that $(\pi, V) \in ob({\mathfrak{R}}^s(G))\Longleftrightarrow$  every irreducible sub-quotient of $\pi$ has inertial support $s=[L, \sigma]_G$.\par

	We have the following theorem, which is Theorem 2.10 in \cite{Bernstein1984}:
	
	\begin{theorem}\label{The_1}
		The Bernstein decomposition gives a direct product decomposition of
		$\mathfrak{R}(G)$ into indecomposable subcategories $\mathfrak{R}^s(G)$:
		\[
		\mathfrak{R}(G) = \prod_{s \in B(G)} \mathfrak{R}^s(G),
		\]
		where $s$ runs over the spectrum $B(G)$. Concretely, if $(\pi, V) \in \mathfrak{R}(G)$,
		then for each $s \in B(G)$, $V$ has a unique maximal $G$-subspace
		$V^s \in \mathfrak{R}^s(G)$ such that
		\[
		V = \bigoplus_{s \in B(G)} V^s.
		\]
		Moreover, if $s, s' \in B(G)$ with $s \neq s'$, then
		\[
		\mathrm{Hom}_G(\mathfrak{R}^s(G), \mathfrak{R}^{s'}(G)) = 0.
		\]
	\end{theorem}
	\section{Bernstein decomposition of Hecke Algebra $\mathcal{H}_G(F)$}
	
	For each equivalence class $s \in \mathcal{B}(G)$, let $\mathcal{H}_G^s(F)$ be the two-sided ideal of $\mathcal{H}_G(F)$ corresponding to all smooth representations $(\pi, V)$ of $G$ with inertial support $s$. In particular, $\mathcal{H}_G^s(F)$ is the unique $G$-subspace of $\mathcal{H}_G(F)$ lying in $\mathcal{R}^s(G)$ such that it is maximal for this property. Each ideal $\mathcal{H}_G^s(F)$ is a non-commutative, non-unital, non-finitely-generated, non-reduced $\mathbb{C}$-algebra and it is characterized by:
	
	\[
	\mathcal{H}_G^s(F)\cdot V = V \quad \text{for} \quad V \in \mathcal{R}^s(G)
	\]
	
	The Bernstein decomposition (Theorem~\ref{The_1}) induces a factorization of the Hecke algebra $\mathcal{H}(G)$ into two-sided ideals as:
	
	\begin{equation}
		\mathcal{H}_G(F) = \bigoplus_{s \in \mathcal{B}(G)} \mathcal{H}_G^s(F).
	\end{equation}
	
	The $\mathbb{C}$-algebra $\mathcal{H}_G^s(F)$ is called the Bernstein block of the Hecke algebra $\mathcal{H}_G(F)$ corresponding to the inertial class $s$. \par
	
	If $\mathcal{H}_G^s(F)\text{-Mod}$ denotes the category of all unitary $\mathcal{H}_G^s(F)$-modules, then there is a natural equivalence of categories:
	
	\begin{equation}\label{eqn_1}
		\mathcal{R}^s(G) \simeq \mathcal{H}_G^s(F)\text{-Mod}.
	\end{equation} 
	
	\section{Types and Spherical Hecke algebras}
	
\subsection{Types}

Let $K$ be a compact open subgroup of $G$. Let $(\rho, W)$ be an irreducible smooth representation of $K$ and $(\pi, V)$  be a  smooth representation of $G$. Let $V^{\rho}$ be the $\rho$-isotopic subspace of $V$. Thus $V^{\rho}$ is the sum of all irreducible $K$-subspaces of $V$ which are equivalent to $\rho$.

\[
V^{\rho}= \sum\limits_{W'}W'
\] where the sum is over all $W'$ such that $(\pi|_K, W') \simeq (\rho, W)$.\par

Recall that $\mathcal{H}_G(F)$ is the space of locally constant compactly supported functions with values in $\mathbb{C}$. Let $e_\rho$ be the element in $\mathcal{H}_G(F)$ with support $K$ such that

\begin{center}
	$e_\rho(x)= \frac{\text{dim}\rho}{\mu(K)}tr_W(\rho(x^{-1})), x \in K.$\\
\end{center}

Note that $e_\rho * e_\rho = e_\rho $ and $e_\rho V = V^{\rho}$. Let $\mathfrak{R}_\rho(G)$ be the full sub-category of $\mathfrak{R}(G)$ consisting of all $(\pi, V)$ where $V$ is generated by the elements of $V^\rho$. So $(\pi, V) \in \mathfrak{R}_\rho(G) \Longleftrightarrow V = \mathcal{H}_G(F) * e_\rho V$. We now state the definition of a $s$-type for $ s \in \mathfrak{B}(G)$.\par

\begin{definition}
	Let $ s \in \mathfrak{B}(G)$. We say that $(K, \rho)$ is an $s$-type in $G$ if $\mathfrak{R}_\rho(G)=
	\mathfrak{R}^s(G)$.
\end{definition}

\subsection{Spherical Hecke algebras}

Recall that $G=\mathrm{GL}_n(F$) and $K$ is a compact open subgroup of $G$. Let $(\rho, W)$ be an irreducible smooth representation of $K$. Here we introduce the spherical Hecke algebra $\mathcal{H}(G,\rho)$.\par

\[
\mathcal{H}(G,\rho)= \left\lbrace f \colon G \to End_{\mathbb{C}}(\rho^{\vee}) \; \middle|  \;
\begin{varwidth}{\linewidth}
	supp($f$) is compact and \\
	$f(k_1gk_2)= \rho^{\vee}(k_1)f(g)\rho^{\vee}(k_2)$\\
	where $k_1,k_2 \in K, g \in G$
\end{varwidth}
\right \rbrace.
\]\par

Then $\mathcal{H}(G,\rho)$ is a $\mathbb{C}$-algebra with multiplication given by convolution $*$ w.r.t some fixed Haar measure $\mu$ on $G$. So for elements $f,g \in \mathcal{H}(G)$ we have 
\[(f * g)(x)= \int_G f(y)g(y^{-1}x)d\mu(y).\] \par

The importance of types is seen from the following result. Let $\pi$ be a smooth representation in $\mathfrak{R}^{s}(G)$. Let $\mathcal{H}(G, \rho)$- Mod denote the category of $\mathcal{H}(G, \rho)$-modules. If $(K, \rho)$ is an $s$-type then 

\begin{equation}\label{eqn_2}
m_G \colon \mathfrak{R}^{s}(G) \longrightarrow \mathcal{H}(G, \rho)\text{-Mod} 
\end{equation}

given by $m_G(\pi)= \mathrm{Hom}_K(\rho,\pi)$ is an equivalence of categories.\par 

From equations (\ref{eqn_1}) and (\ref{eqn_2}), it follows that the given an $s \in B(G)$, there exists $(K, \rho)$ of $s$-type such that the following equivalence of categories hold good:

\begin{equation}\label{eqn_3}
\mathcal{H}^s_G(F)\text{-Mod}	\xrightarrow{\sim} \mathcal{H}(G, \rho)\text{-Mod}
\end{equation}

By the Main Theorem of \cite{BushnellKutzko1999}, the Hecke algebra $\mathcal{H}(G, \rho)$ is isomorphic to the tensor product\[\bigotimes_{i=1}^{r} \mathcal{H}(n_i,q^{k_i})\]of affine Hecke algebras. Here, $q$ is the cardinality of the residue field of \(K\), while \(k_i\) denotes the torsion number  of \(\rho\) (where $\rho= \rho_1 \otimes \cdots \otimes \rho_r$), cf.~\cite[p.~22]{AubertBaumPlymen}. In particular,\[q^{k_i} \neq -1\]for all \(i\).

From \ref{eqn_3} and the above discussion we have,

\begin{equation}\label{eqn_4}
	\mathcal{H}^s_G(F)\text{-Mod}	\xrightarrow{\sim} \mathcal{H}(G, \rho)\text{-Mod} \xrightarrow{\sim} \bigotimes_{i=1}^{r} \mathcal{H}(n_i,q^{k_i})\text{-Mod}
\end{equation}

\section{Final Result}
Let us recall affine Hecke algebras. Let $\mathcal{H}(m,l)$ and $\mathcal{H}(m',l')$ be two affine Hecke Algebras where $m,m' \in \mathbb{Z_+}$ and $r,r' \in \mathbb{C}^{\times}$. Let the generators of $\mathcal{H}(m,l)$ be $S_i(1 \leq i \leq m), T, T^{-1}$ and the generators of $\mathcal{H}(m',l')$ be $s_i(1 \leq i \leq m), t, t^{-1}$ In \cite{Repaka2026}, on page-4 we have shown that if $\mathcal{H}(m,l) \cong \mathcal{H}(m',l')$ as $\mathbb{C}$ algebras then $m=m'$ and $l=l'$. Namely, if $\varphi: \mathcal{H}(m,l) \rightarrow \mathcal{H}(m',l')$ is a $\mathbb{C}$-algebra isomorphism, then it turns out $\varphi(S_i) = s_j$ for some  $j$ such that $1 \leq j \leq m$ and $\varphi$ is a $L_1$-isometry and $m=m', l=l'$. \par

Note that corresponding equation for non-Archimedean local field $E$ as that of equation (\ref{eqn_4}) would be:

\begin{equation}\label{eqn_5}
	\mathcal{H}^s_G(E)\text{-Mod}	\xrightarrow{\sim} \mathcal{H}(G, \rho)\text{-Mod} \xrightarrow{\sim} \bigotimes_{i=1}^{r'} \mathcal{H}(n'_i,q^{k'_i})\text{-Mod}
\end{equation}

As $\mathcal{H}_G(F)-\text{Mod} \cong \mathcal{H}_G(E)-\text{Mod}$, hence from equation (\ref{eqn_4}) and equation (\ref*{eqn_5}) we have,

\begin{equation}
\bigotimes_{i=1}^{r} \mathcal{H}(n_i,q^{k_i})\text{-Mod} \xrightarrow{\sim} \bigotimes_{i=1}^{r'}\mathcal{H}(n'_i,q^{k'_i})
\text{-Mod}
\end{equation}

The above equation is possible if and only if $r=r', n_i=n_i'$ for all $i$ and $q^{k_i}=q^{k'_i}$ for all $i$. Hence, $\bigotimes_{i=1}^{r}\mathcal{H}(n_i,q^{k_i}) \sim  \bigotimes_{i=1}^{r'}\mathcal{H}(n'_i,q^{k'_i})$ which would imply that $\mathcal{H}(G, \rho)$-Mod does not depend on underlying field $F$ and hence $\mathcal{H}_G$-Mod does not depend on underlying field $F$ which proves Theorem \ref{The_0}.

\end{document}